\documentclass[titlepage,12pt]{article}
\usepackage{amssymb}
\usepackage{amsfonts}
\usepackage{amsmath}
\begin{document}

{\Large \bf A Fixed Point Conjecture} \\ \\

{\bf Elem\'{e}r E Rosinger} \\
Department of Mathematics \\
and Applied Mathematics \\
University of Pretoria \\
Pretoria \\
0002 South Africa \\
eerosinger@hotmail.com \\ \\

{\bf Abstract} \\

Inverse limits, unlike direct limits, can in general be void, [1]. The existence of fixed
points for arbitrary mappings $T : X \longrightarrow X$ is conjectured to be equivalent with
the fact that related direct limits of all finite partitions of X are not void. \\ \\

{\bf 1. The Setup} \\ \\

Let $X$ be a nonvoid set and $T : X \longrightarrow X$ a mapping. We denote by \\

(1)~~~ $ {\cal FP} ( X ) $ \\

the set of all finite partitions of $X$. \\

Given $x \in X$ and $\Delta \in {\cal FP} ( X )$, then obviously \\

(2)~~~ $ \exists~~ A \in \Delta ~:~
                  \{~ n \in \mathbb{N}_+ ~|~ T^n ( x ) \in A ~\} ~~\mbox{is infinite} $ \\

Here and in the sequel, we use the notation $\mathbb{N} = \{ 0, 1, 2, 3, \ldots \}$ and
$\mathbb{N}_+ = \{ 1, 2, 3, \ldots \}$. \\

Let us therefore denote \\

(3)~~~ $ \Delta( x ) = \{~ A \in \Delta ~|~
                  \{~ n \in \mathbb{N}_+ ~|~ T^n ( x ) \in A ~\} ~~\mbox{is infinite} ~\} $ \\

In view of (2) we obtain \\

(4)~~~ $ \Delta( x ) \neq \phi $ \\

In setting up the fixed point conjecture, it is useful to consider the following two simple
instances. \\

{\bf Example 1} \\

1) Let $T = id_X$, that is, the {\it identity} mapping on $X$. Then for $x \in X$ and $A
\subseteq X$, we clearly have \\

(5)~~~ $ \{~ n \in \mathbb{N}_+ ~|~ T^n ( x ) \in A ~\} ~=~
                             \begin{array}{|l}
                                 ~\phi ~~\mbox{if}~ x \notin A \\ \\
                                 ~\mathbb{N}_+ ~~\mbox{if}~ x \in A
                              \end{array} $ \\

hence for $\Delta \in {\cal FP} ( X )$, we obtain \\

(6)~~~ $ \Delta ( x ) ~=~ \{~ A ~\}, ~~\mbox{where}~ x \in A $ \\

2) Let $T$ be a {\it constant} mapping on $X$, that is, $T ( x ) = c$, for $x \in X$, where
$c \in X$ is given. Then for $x \in X$ and $A \subseteq X$, we clearly have \\

(7)~~~ $ \{~ n \in \mathbb{N}_+ ~|~ T^n ( x ) \in A ~\} ~=~
                             \begin{array}{|l}
                                 ~\phi ~~\mbox{if}~ c \notin A \\ \\
                                 ~\mathbb{N}_+ ~~\mbox{if}~ c \in A
                              \end{array} $ \\

hence for $\Delta \in {\cal FP} ( X )$, we obtain \\

(8)~~~ $ \Delta ( x ) ~=~ \{~ A ~\}, ~~\mbox{where}~ c \in A $

\hfill $\Box$ \\

Let us recall now the following natural partial order structure on ${\cal FP} ( X )$ given by
the concept of {\it refinement} of partitions. Namely, if $\Delta , \Delta\,' \in {\cal FP}
( X )$, then we denote \\

(9)~~~ $ \Delta \leq \Delta\,' $ \\

if and only if \\

(10)~~~ $ \forall~~ A\,' \in \Delta\,' ~:~
                        \exists~~ A \in \Delta ~:~ A\,' \subseteq A $ \\

and in such a case, we define the mapping \\

(11)~~~ $ \psi_{\Delta\,',\, \Delta} : \Delta\,' \longrightarrow \Delta\ $ \\

by \\

(12)~~~ $ A\,' \subseteq A = \psi_{\Delta\,',\, \Delta} ( A\,' ) $ \\

Then obviously \\

(13)~~~ $ \psi_{\Delta\,',\, \Delta} ( \Delta\,'( x ) ) \subseteq \Delta( x ) $ \\

Indeed, if $A\,' \in \Delta\,'( x )$, then (3) gives \\

$~~~~~~ \{~ n \in \mathbb{N}_+ ~|~ T^n ( x ) \in A\,' ~\} ~~\mbox{is infinite} $ \\

but in view of (10), we have \\

$~~~~~~ A\,' \subseteq \psi_{\Delta\,',\, \Delta} ( A\,' ) $ \\

hence \\

$~~~~~~ \{~ n \in \mathbb{N}_+ ~|~ T^n ( x ) \in \psi_{\Delta\,',\, \Delta} ( A\,' ) ~\}
                                                                  ~~\mbox{is infinite} $ \\

thus (13). \\

We note now that the finite partitions \\

(14)~~~ $ \Delta \in {\cal FP} ( X )$ \\

together with the mappings \\

(15)~~~ $ \psi_{\Delta\,',\, \Delta},~~~ \Delta,~ \Delta\,' \in {\cal FP} ( X ),~ \Delta \leq
                                                                 \Delta\,' $ \\

form an {\it inverse} family, [1, p.191]. \\

Furthermore, in view of (13), we also have the following stronger version of the above. For
every $x \in X$, let us define the mappings \\

(16)~~~ $ \psi_{\Delta\,',\, \Delta,\, x} : \Delta\,' ( x ) \longrightarrow \Delta ( x ),~~~
                          \Delta,~ \Delta\,' \in {\cal FP} ( X ),~ \Delta \leq \Delta\,' $ \\

by \\

(17)~~~ $ \psi_{\Delta\,',\, \Delta,\, x} = \psi_{\Delta\,',\, \Delta} |_{\Delta\,' ( x )} $ \\

Then again, for every $x \in X$, we obtain the {\it inverse} family \\

(18)~~~  $ \Delta \in {\cal FP} ( X )$ \\

(19)~~~ $ \psi_{\Delta\,',\, \Delta,\, x} : \Delta\,' ( x ) \longrightarrow \Delta ( x ),~~~
                          \Delta,~ \Delta\,' \in {\cal FP} ( X ),~ \Delta \leq \Delta\,' $ \\

Consequently, for each $x \in X$, we can consider the {\it inverse limit} \\

(20)~~~ $ \underleftarrow{\lim}_{\,\Delta \,\in\, {\cal FP} ( X )}~ \Delta ( x ) $ \\

Returning now to the two simple instances in Example 1 above, we further have \\

{\bf Example 2} \\

1) In the case 1) of Example 1, it follows easily that, for $x \in X$, we have \\

(21)~~~ $ \underleftarrow{\lim}_{\,\Delta \,\in\, {\cal FP} ( X )}~ \Delta ( x ) = \{~ ( \xi_\Delta ~|~ \Delta \in {\cal FP} ( X ) ) ~\} $ \\

where \\

(22)~~~ $ \xi_\Delta = x,~~~ \Delta \in {\cal FP} ( X ) $ \\

2) In the case 2) of Example 1, for $x \in X$, we easily obtain \\

(23)~~~ $ \underleftarrow{\lim}_{\,\Delta \,\in\, {\cal FP} ( X )}~ \Delta ( x ) ~=~
                             \begin{array}{|l}
                                 ~\phi ~~\mbox{if}~ x \neq c  \\ \\
                                 ~\{~ ( \xi_\Delta ~|~ \Delta \in {\cal FP} ( X ) ) ~\}
                                                                      ~~\mbox{if}~ x = c
                              \end{array} $ \\

where \\

(24)~~~ $ \xi_\Delta = c,~~~ \Delta \in {\cal FP} ( X ) $ \\

{\bf Remark 1} \\

At this stage, an important fact to note is that, in general, an inverse limit such as in (20)
may be void, [1, (c) in Exercise 4, p. 252], even if none of the sets $\Delta ( x )$ is void,
and each the mappings $\psi_{\Delta\,',\, \Delta,\, x}$ is surjective. Therefore, the
relations (21) and (23), even if easy to establish, are as inverse limits nontrivial, in view of the
arbitrariness of the sets $X$ and mappings $T$ involved. \\

A common feature of the mappings $T : X \longrightarrow X$ in both cases above is that they
have {\it fixed points}. Namely, for the identity mapping $T = id_X$, each point $x \in X$ is
such a fixed point, while for the constant mapping $T = c$, the point $x = c \in X$ is the
only fixed point. \\

Further, as suggested by Scott Kominers, Daniel Litt and Brett Harrison, in view of the fact
that a fixed point of a mapping $T : X \longrightarrow X$ is but a {\it particular} case of a
periodic point of that mapping, or equivalently, of a fixed point of the mapping $T^n$, for
some $n \in \mathbb{N}_+$, we are led to the \\

{\bf Conjecture} \\

Given a nonvoid set $X$ and a mapping $T : X \longrightarrow X$, then for $x \in X$, we have \\

(28)~~~ $ \underleftarrow{\lim}_{\,\Delta \,\in\, {\cal FP} ( X )}~ \Delta ( x ) \neq \phi
                           ~~~\Longleftrightarrow~~~
                                  (~~ \exists~~~ n \in \mathbb{N}_+ ~:~ T^n ( x ) = x ~~) $ \\

where the issue is whether the implication $\Longrightarrow$ holds, since the converse
implication is easy to establish.

\end{document}